\documentclass[11pt,a4paper,fullpage]{article}

\usepackage{amsmath,amsfonts,amssymb,amsthm,graphics,amscd}
\usepackage{dsfont}
\newtheorem{Teorema}{Theorem}[section]
\newtheorem{Propozicija}[Teorema]{Proposition}
\newtheorem{Definicija}[Teorema]{Definition}
\newtheorem{Posledica}[Teorema]{Corollary}

\newtheorem{Lema}[Teorema]{Lemma}
\newtheorem{Primedba}[Teorema]{Remark}

\numberwithin{equation}{section}

\begin{document}
	\title {Invertible completions of FLI and FRI upper triangular operator matrices}

\author{Nikola Sarajlija\footnote{corresponding author: Nikola Sarajlija; University of Novi Sad, Faculty of Sciences, Novi Sad 21000, Serbia; {\it e-mail}: {\tt nikola.sarajlija@dmi.uns.ac.rs}} \footnote{Author is supported by the Ministry of Education, Science and Technological Development of the Republic of Serbia under grants no. 451-03-33/2026-03/ 200125 and 451-03-34/2026-03/ 200125}
}
\maketitle

\begin{abstract}
If $A\in\mathcal{B}(\mathcal{H})$ and $B\in\mathcal{B}(\mathcal{K})$ are given operators, denote by $M_C$ an operator matrix of the form $$M_C=\begin{pmatrix}
    A & C\\ 0 & B
\end{pmatrix}\in\mathcal{B}(\mathcal{H}\oplus\mathcal{K})$$ acting on a direct sum of infinite dimensional separable Hilbert spaces $\mathcal{H}$ and $\mathcal{K}$, where $C\in\mathcal{B}(\mathcal{K},\mathcal{H})$ is unknown. In this article we solve the completion problem of $M_C$ to Fredholm left and Fredholm right invertibility, and we obtain appropriate perturbation results as consequences. We illustrate our results by solving the 'filling in holes' problem for Fredholm left and Fredholm right spectra. Finally, we consider some special classes of operators. Throughout the article, we recover some known results.
\end{abstract}

\textit{$2020$ Math. Subj. Class:} 47A08, 47A05, 47A53, 47A55

\vspace{2mm}

\section{Introduction and Preliminaries}

Throughout this text, let $\mathcal{H}$ and $\mathcal{K}$ be infinite dimensional separable Hilbert spaces. Let $\mathcal{B}(\mathcal{H},\mathcal{K})$ stand for the collection of all bounded linear operators mapping $\mathcal{H}$ to $\mathcal{K}$. We write $\mathcal{B}(\mathcal{H})=\mathcal{B}(\mathcal{H},\mathcal{H})$ for short. If $T\in\mathcal{B}(\mathcal{H},\mathcal{K})$, then $\mathcal{N}(T)$ and $\mathcal{R}(T)$ stand for the null space and the range space of $T$, respectively. Denote by $\alpha(T)=\dim\mathcal{N}(T)$ and $\beta(T)=\dim Y/\mathcal{R}(T)$ corresponding Hilbert space dimensions. Notions $\alpha(\cdot)$ and $\beta(\cdot)$ are called nullity and deficiency, respectively, and are used extensively in Fredholm theory. We state a well known result that we use in the sequel.
\begin{Lema}\label{ZATVORENALEMA}
If $T\in\mathcal{B}(\mathcal{H},\mathcal{K})$, then the following implication holds:
$$\beta(T)<\infty \Rightarrow\ \mathcal{R}(T)\ is\ closed.$$
\end{Lema}
We use the standard notation $T^*$ to denote the adjoint operator of $T$. Operators $T$ and $T^*$ show dual behaviour in terms of spectral properties, as the following simple results shows:
\begin{Lema}\label{DUALNALEMA}
If $T\in\mathcal{B}(\mathcal{H},\mathcal{K})$ then the following holds:
\begin{itemize}
    \item[(a)] $\alpha(T)=\beta(T^*)$;
    \item[(b)] $\beta(T)=\alpha(T^*)$.
\end{itemize}
\end{Lema}
We will use Lemma \ref{DUALNALEMA} in duality arguments.\\

Next, we say that $T$ is an upper (lower) semi-Fredholm operator if $\alpha(T)<\infty$  ($\beta(T)<\infty$) and $\mathcal{R}(T)$ is closed in $\mathcal{K}$. If $T$ is an upper (lower) semi-Fredholm operator with $\alpha(T)=0$ ($\beta(T)=0$), then $T$ is a left 
 (right) invertible operator. Next, $T$ is a Fredholm operator if $\alpha(T),\beta(T)<\infty$, and $T$ is invertible if $T$ is Fredholm invertible with $\alpha(T)=\beta(T)=0$. \\

Corresponding spectra of an operator $T\in\mathcal{B}(\mathcal{H},\mathcal{K})$ are defined as follows: the left spectrum is $\sigma_{l}(T)=\lbrace\lambda\in\mathds{C}: \lambda-T\ is\ not\ left\ invertible\rbrace$, the right spectrum is $\sigma_{r}(T)=\lbrace\lambda\in\mathds{C}: \lambda-T\ is\ not\ right\ invertible\rbrace$,
while the spectrum of $T$ is set $\sigma(T)=\lbrace\lambda\in\mathds{C}: \lambda-T\ is\ not\ invertible\rbrace$. All of these spectra are compact nonempty subsets of the complex plane. We write $\rho_{l}(T), \rho_{r}(T), \rho(T)$ for their complements, respectively. Next, we use the following notation for the upper and lower semi-Fredholm spectrum, respectively: $\sigma_{SF+}(T)=\lbrace\lambda\in\mathds{C}:\ \lambda-T\ is\ not\ upper\ semi-Fredholm\rbrace$, $\sigma_{SF-}(T)=\lbrace\lambda\in\mathds{C}:\ \lambda-T\ is\ not\ lower\ semi-Fredholm \rbrace.$ Appropriate resolvent sets are denoted by $\rho_{SF+}(T)$ and $\rho_{SF-}(T)$, respectively. Finaly, $\sigma_p(T)=\lbrace\lambda\in\mathds{C}:\ \lambda-T\ is\ not\ injective\rbrace$ is the point spectrum of $T$, and $\sigma_d(T)=\lbrace\lambda\in\mathds{C}:\ \lambda-T\ is\ not\ surjective\rbrace$ is the defect spectrum of $T$. Notice that $\sigma(T)=\sigma_l(T)\cup\sigma_r(T)=\sigma_p(T)\cup\sigma_d(T)$.\\

All of these types of invertibility are introduced some time ago and are well studied. However, in this article we investigate a different type of invertibility that was introduced very recently in \cite{DONG2}. Therefore, we state the following definition.

\begin{Definicija} \cite{DONG2}
    Let $T\in\mathcal{B}(\mathcal{H},\mathcal{K})$. Then $T$ is a Fredholm left invertible operator (FLI operator for short) if $\alpha(T)=0$ and $\beta(T)<\infty$.
\end{Definicija}
Dually, we introduce here the following definition:
\begin{Definicija}
    Let $T\in\mathcal{B}(\mathcal{H},\mathcal{K})$. Then $T$ is a Fredholm right invertible operator (FRI operator for short) if $\alpha(T)<\infty$ and $\beta(T)=0$.
\end{Definicija}
We denote the corresponding spectra and resolvent sets by $\sigma_{FLI}(T),\sigma_{FRI}(T),\rho_{FLI}(T),\rho_{FRI}(T)$.\\

One way of studying spectral properties of an operator acting on a Hilbert space $\mathcal{H}$, is to consider its restrictions to certain closed subspaces of $\mathcal{H}$, and to analyze operators obtained in that manner. This led to the notion of operator matrices, which have been widely investigated from various points of view thus far (see \cite{DUNFORD}, \cite{TRETTER}). In this article, we focus on a special type of operator matrices, namely $2\times2$ upper triangular operator matrices.   \\

Therefore, if $A\in\mathcal{B}(\mathcal{H})$ and $B\in\mathcal{B}(\mathcal{K})$ are given operators, denote by $M_C$ an operator matrix of the form $$M_C=\begin{pmatrix}
    A & C\\ 0 & B
\end{pmatrix}\in\mathcal{B}(\mathcal{H}\oplus\mathcal{K})$$ where $C\in\mathcal{B}(\mathcal{K},\mathcal{H})$ is unknown. In this article we find necessary and sufficient conditions for $A$ and $B$ so that exists $C\in\mathcal{B}(\mathcal{K},\mathcal{H})$ such that $M_C$ is FLI or FRI, thus solving two completion problems for $M_C$. Spectral properties and completion problems of $M_C$ have been studied well in the past (see \cite{CAO}, \cite{OPERATORTHEORY}, \cite{DU}, \cite{WANG}, etc.), and very recently there have been some extensions to upper triangular operator matrices of dimension 3 (\cite{DONG},\cite{SARAJLIJA4}) and of general dimension $n>3$ \cite{SARAJLIJA2}. \\

Before proceeding any further, we state two preliminary results.
\begin{Lema}\label{JEDNOSTAVNALEMA}
Let $S,T\in\mathcal{B}(\mathcal{H},\mathcal{K})$. If $T$ is invertible, then:
\begin{itemize}
\item[(a)] $\mathcal{R}(TS)\cong\mathcal{R}(S)$;
\item[(b)] $\mathcal{N}(TS)\cong\mathcal{N}(S)$.
\end{itemize}
\end{Lema}

\begin{Lema}\label{POMOCNALEMA}
In the notation as above, the following is true:
\begin{itemize}
\item[(a)] If $M_C$ is left invertible, then $A$ is left invertible;
\item[(b)]  If $M_C$ is right invertible, then $B$ is right invertible;
\item[(c)]  If $M_C$ is upper semi-Fredholm, then $A$ is upper semi-Fredholm;
\item[(d)] If $M_C$ is lower semi-Fredholm, then $B$ is lower semi-Fredholm;
\item[(e)] If $M_C$ is Fredholm, then $A$ is Fredholm if and only if $B$ is Fredholm.
\end{itemize}
\end{Lema}
\textbf{Proof.} All statements easily follow by definitions introduced previously, and with regards to matrix decomposition
\begin{equation}\label{RAZLAGANJE}
M_C=\begin{pmatrix}I_{\mathcal{H}} & 0\\0 & B\end{pmatrix}\begin{pmatrix}I_{\mathcal{H}} & C\\0 & I_{\mathcal{K}}\end{pmatrix}\begin{pmatrix}A & 0\\0 & I_{\mathcal{K}}\end{pmatrix}.\quad \square
\end{equation}\\

Article is organized as follows. In this first section we have explained the topic of our investigation and we have provided some preliminary results. In the second section we successfully solve completion problems to Fredholm left and Fredholm right invertibility, and obtain perturbation results as consequences. Finally, in the third section, we solve the 'filling in holes' problem for Fredholm left and Fredholm right spectra, and we consider some special classes of operators.

\section{Completion problems of $M_{C}$ to Fredholm left and Fredholm right invertibility}

In this section we provide necessary and sufficient conditions for Fredholm left and Fredholm right invertibility of $M_{C}$. 

\begin{Teorema}\label{GLAVNA1}
Let $A\in\mathcal{B}(\mathcal{H}),B\in\mathcal{B}(\mathcal{K})$ be given. There exists $C\in\mathcal{B}(\mathcal{K},\mathcal{H})$ such that $M_C$ is a FLI operator if and only if the following holds: 
\begin{itemize}
\item[(a)] $A$ is left invertible;
\item[(b)]  $B$ is lower semi-Fredholm;
\item[(c)]  $\alpha(B)\leq\beta(A)<\infty$ or $\alpha(B)=\beta(A)=\infty$.
\end{itemize}
\end{Teorema}
\textbf{Proof.} Assume that $M_C$ is a FLI operator. Then, $M_C$ is left invertible, hence part $(a)$ follows by Lemma \ref{POMOCNALEMA}$(a)$. Furthermore, $M_C$ is lower semi-Fredholm, hence part $(b)$ follows by Lemma \ref{POMOCNALEMA}$(d)$. Let us prove part $(c)$. Consider upper triangular operators $\begin{pmatrix}I_\mathcal{H} & 0\\ 0 & B\end{pmatrix}$ and $\begin{pmatrix}I_{\mathcal{H}} & C\\ 0 & I_{\mathcal{K}}\end{pmatrix}\begin{pmatrix}A & 0\\0 & I_{\mathcal{K}}\end{pmatrix}$ whose product is equal to $M_C$. Lemma \ref{ZATVORENALEMA} and assumption $(a)$ guarantee that their range is closed, so one can apply Harte's index theorem \cite{HARTE2} to obtain
\begin{equation}
\begin{aligned}
    \mathcal{N}\Bigg(\begin{pmatrix}I & 0\\ 0 & B\end{pmatrix}\Bigg)\times\mathcal{N}\Bigg(\begin{pmatrix}I & C\\ 0 & I\end{pmatrix}\begin{pmatrix}A & 0\\0 & I\end{pmatrix}\Bigg)\times\mathcal{R}(M_C)^\perp\cong \\ \mathcal{N}(M_C)\times\mathcal{R}\Bigg(\begin{pmatrix}I & 0\\ 0 & B\end{pmatrix}\Bigg)^\perp\times\mathcal{R}\Bigg(\begin{pmatrix}I & C\\ 0 & I\end{pmatrix}\begin{pmatrix}A & 0\\0 & I\end{pmatrix}\Bigg)^\perp.
\end{aligned}
\end{equation}
This further implies
\begin{equation}
    \mathcal{N}(B)\times\mathcal{N}(A)\times\mathcal{R}(M_C)^\perp\cong\mathcal{N}(M_C)\times\mathcal{R}(B)^\perp\times\mathcal{R}(A)^\perp,
\end{equation}
where we used Lemma \ref{JEDNOSTAVNALEMA}. 
Left invertibility of $M_C$ and left invertibility of $A$ imply
\begin{equation}
    \mathcal{N}(B)\times\mathcal{R}(M_C)^\perp\cong\mathcal{R}(B)^\perp\times\mathcal{R}(A)^\perp.
\end{equation}
which finally yields $\alpha(B)+\beta(M_C)=\beta(B)+\beta(A)$. Since $M_C$ is FLI one can deduce that $\infty>\beta(M_C)\geq\beta\begin{pmatrix}I & 0\\ 0 & B\end{pmatrix}=\beta(B)$. Therefore, either $\alpha(B)=\beta(A)=\infty$ or $\alpha(B)\leq\beta(A)<\infty$. \\

To prove the converse implication, assume that parts $(a)-$ $(c)$ are valid. Assumption $(a)$ implies that $\mathcal{R}(A)$ is closed, hence $\mathcal{H}=\mathcal{N}(B)\oplus\mathcal{N}(B)^\perp$ and $\mathcal{K}=\mathcal{R}(A)^\perp\oplus\mathcal{R}(A)$. \\Assume first that $\alpha(B)\leq\beta(A)<\infty$. If $\alpha(B)=0$ we choose $C=\mathbf{0}$, therefore assume $0<\alpha(B)\leq\beta(A)<\infty$. We conclude that there exists a left invertible finite rank operator $J:\mathcal{N}(B)\to\mathcal{R}(A)^\perp$. Put $C=\begin{pmatrix}
    J & 0\\ 0 & 0
\end{pmatrix}:\begin{pmatrix}
    \mathcal{N}(B) \\ 
    \mathcal{N}(B)^\perp
\end{pmatrix}\to\begin{pmatrix}
    \mathcal{R}(A)^\perp\\\mathcal{R}(A)
\end{pmatrix}$. We shall prove that such choice of $C$ gives an FLI operator $M_C$. First, we prove that $M_C$ is injective. Put $M_C\begin{pmatrix} x \\ y\end{pmatrix}=\begin{pmatrix}0 \\ 0\end{pmatrix}$. Then $By=0$, which implies $y\in\mathcal{N}(B)$. Also, $Ax+Cy=0$, which together with $y\in\mathcal{N}(B)$ implies $Ax+Jy=0$. Since $\mathcal{R}(A)\cap\mathcal{R}(J)=\lbrace0\rbrace$ we have $Ax=Jy=0$, and finally $x=y=0$. We have thus proved that $\alpha(M_C)=0$. Now, we shall prove that $\beta(M_C)<\infty$. Notice that $\mathcal{R}(M_C)=(\mathcal{R}(A)\oplus\mathcal{R}(J))\oplus\mathcal{R}(B)$. Therefore, $(\mathcal{H}\oplus\mathcal{K})/{\mathcal{R}(M_C)\cong\mathcal{R}(J)^\perp\oplus\mathcal{R}(B)^\perp}$, which implies $\beta(M_C)\leq\beta(A)+\beta(B)<\infty$.\\ Assume next $\alpha(B)=\beta(A)=\infty$. Then there is an isomorphism $J:\mathcal{N}(B)\to\mathcal{R}(A)^\perp$. In this case we choose $M_C$ analogously as in the previous case. Now, injectivity of $M_C$ is proved in the same fashion, and $\mathcal{R}(M_C)=X\oplus\mathcal{R}(B)$ implying $\beta(M_C)=\beta(B)<\infty.$ $\square$\\

By duality, it is easy to obtain an analogous result for the Fredholm right invertibility of $M_C$.

\begin{Teorema}\label{GLAVNA2}
Let $A\in\mathcal{B}(\mathcal{H}),B\in\mathcal{B}(\mathcal{K})$ be given. There exists $C\in\mathcal{B}(\mathcal{K},\mathcal{H})$ such that $M_C$ is a FRI operator if and only if the following holds: 
\begin{itemize}
\item[(a)] $B$ is right invertible;
\item[(b)]  $A$ is upper semi-Fredholm;
\item[(c)]  $\beta(A)\leq\alpha(B)<\infty$ or $\alpha(B)=\beta(A)=\infty$.
\end{itemize}
\end{Teorema}
\textbf{Proof. }This follows directly from Theorem \ref{GLAVNA1} and Lemma \ref{DUALNALEMA}. $\square$\\

An obvious consequence of Theorem \ref{GLAVNA1} and Theorem \ref{GLAVNA2} is the following well known result.
\begin{Propozicija}\cite{DU}
Let $A\in\mathcal{B}(\mathcal{H}),B\in\mathcal{B}(\mathcal{K})$ be given. There exists $C\in\mathcal{B}(\mathcal{K},\mathcal{H})$ such that $M_C$ is an invertible operator if and only if the following holds: 
\begin{itemize}
\item[(a)] $A$ is left invertible;
\item[(b)] $B$ is right invertible;
\item[(c)] $\alpha(B)=\beta(A)$.
\end{itemize}
\end{Propozicija}
Now, we are able to pose appropriate perturbation results.

\begin{Posledica}\label{POSLEDICA10'}
Let $A\in\mathcal{B}(\mathcal{H}),\ B\in\mathcal{B}(\mathcal{K})$. Then
$$
\bigcap_{C\in\mathcal{B}(\mathcal{K},\mathcal{H})}\sigma_{FLI}(M_C)=\sigma_{l}(A)\cup\sigma_{SF-}(B)\cup\Delta_{FLI},
$$
where
\begin{equation*}
\begin{aligned}
\Delta_{FLI}=\Big\lbrace\lambda\in\mathds{C}:\ \alpha(B-\lambda)>\beta(A-\lambda)\ or\ \beta(A-\lambda)=\infty
\Big\rbrace\bigcap\\ \Big\lbrace\lambda\in\mathds{C}:\ \alpha(B-\lambda)\neq\beta(A-\lambda)\ or\ \alpha(B-\lambda)=\beta(A-\lambda)<\infty\Big\rbrace.
\end{aligned}
\end{equation*}
\end{Posledica}
Dually, we have:
\begin{Posledica}\label{POSLEDICA10''}
Let $A\in\mathcal{B}(\mathcal{H}),\ B\in\mathcal{B}(\mathcal{K})$. Then
$$
\bigcap_{C\in\mathcal{B}(\mathcal{K},\mathcal{H})}\sigma_{FRI}(M_C)=\sigma_{r}(B)\cup\sigma_{SF+}(A)\cup\Delta_{FRI},
$$
where
\begin{equation*}
\begin{aligned}
\Delta_{FRI}=\Big\lbrace\lambda\in\mathds{C}:\ \beta(A-\lambda)>\alpha(B-\lambda)\ or\ \alpha(B-\lambda)=\infty)\Big\rbrace\bigcap\\ \Big\lbrace\lambda\in\mathds{C}:\ \beta(A-\lambda)\neq\alpha(B-\lambda)\ or\ \beta(A-\lambda)=\alpha(A-\lambda)<\infty\Big\rbrace.
\end{aligned}
\end{equation*}
\end{Posledica}

\section{'Filling in holes' results for $\sigma_{FLI}$ and $\sigma_{FRI}$}

We now move our focus to 'filling in holes' results for Fredholm left and Fredholm right spectrum. We mention here that various 'filling in holes'-type results for left (right) spectrum, upper (lower) semi-Fredholm spectrum and upper (lower) semi-Weyl spectrum can be seen in \cite{SARAJLIJA2}. To begin with, we state the following result.
\begin{Lema}\label{IZMEDJUCEGAJESPEKTAR}
    Let $A\in\mathcal{B}(\mathcal{H}),\ B\in\mathcal{B}(\mathcal{K})$. Then
    \begin{equation}
    \sigma_{FLI}(A)\setminus(\sigma_d(A)\cap\sigma_l(B))\subseteq\sigma_{FLI}(M_C)\subseteq\sigma_{FLI}(A)\cup\sigma_{FLI}(B)
    \end{equation}
    for every $C\in\mathcal{B}(\mathcal{K},\mathcal{H})$.
\end{Lema}
\textbf{Proof. }If $A$ and $B$ are Fredholm, then it is well known that $M_C$ is Fredholm regardless of $C$ (see e.g. \cite[Proposition 3.1]{OPERATORTHEORY}). If $A$ and $B$ are both left invertible, then factorization \eqref{RAZLAGANJE} shows that $M_C$ is left invertible as well (the second factor in (\ref{RAZLAGANJE}) is invertible regardless of $C$). Thus, the second inclusion in (\ref{IZMEDJUCEGAJESPEKTAR}) is proved. To prove the first inclusion, assume  $\lambda\not\in\sigma_{FLI}(M_C)$. Let $\lambda\in\sigma_{FLI}(A)\setminus(\sigma_d(A)\cap\sigma_l(B))$. Then $M_C-\lambda$ is FLI operator, implying that $A-\lambda$ is left invertible. But, $\lambda\in\sigma_{FLI}(A)$, therefore $\beta(A-\lambda)=\infty$. This condition, together with $\lambda\not\in\sigma_d(A)\cap\sigma_l(B)$ gives $\lambda\in\rho_l(B)$. However, this means that $M_C-\lambda$ and $B-\lambda$ are both Fredholm, implying that $A-\lambda$ must be Fredholm as well (see Lemma \ref{POMOCNALEMA}(e)). But, this is contrary to $\beta(A-\lambda)=\infty$. Thus, the first inclusion is proved as well. $\square$ 

This result suggests that the passage from $\sigma_{FLI}(M_0)=\sigma_{FLI}(A)\cup\sigma_{FLI}(B)$ to $\sigma_{FLI}(M_C)$ is accomplished by removing certain subsets that are contained in $\sigma_d(A)\cap\sigma_l(B)$ and disjoint from $\sigma_{FLI}(A)$. The next result confirms this hypothesis.

\begin{Teorema}\label{TEOREMAZAFLI}
    Let $A\in\mathcal{B}(\mathcal{H}),\ B\in\mathcal{B}(\mathcal{K})$. Then, for every $C\in\mathcal{B}(\mathcal{K},\mathcal{H})$,
    $$\eta(\sigma_{FLI}(A)\cup\sigma_{FLI}(B))=\eta(\sigma_{FLI}(M_C)),$$
    where $\eta(\cdot)$ denotes the 'polynomially-convex hull'.
    More precisely,
    $$\sigma_{FLI}(A)\cup\sigma_{FLI}(B)=\sigma_{FLI}(M_C)\cup W,$$
    where $W$ is equal to the union of some holes in $\sigma_{FLI}(A)$ which happen to be subsets of $\sigma_d(A)\cap\sigma_l(B)$.
\end{Teorema}

\begin{Primedba}
    By a 'hole' in $K\subseteq\mathds{C}$ we understand a nonempty bounded component of $\mathds{C}\setminus K$.
\end{Primedba}

\textbf{Proof. } The statement will easily follow once we establish an equality
\begin{equation}\label{JEDNAKOSTSAETOM}
    \eta(\sigma(T))=\eta(\sigma_{FLI}(T)).
\end{equation}
However, obviously $int\ \sigma_{FLI}(T)\subseteq int\ \sigma(T)$ and $\partial\sigma(T)\subseteq\sigma_l(T)\subseteq\sigma_{FLI}(T)$, which implies $\partial\sigma(T)\subseteq\partial\sigma_{FLI}(T)$, altogether proving (\ref{JEDNAKOSTSAETOM}). Assume $\lambda\in(\sigma_{FLI}(A)\cup\sigma_{FLI}(B))\setminus\sigma_{FLI}(M_C)$. Then, according to the first inclusion in Lemma \ref{IZMEDJUCEGAJESPEKTAR} we have $\lambda\in\sigma_d(A)\cap\sigma_l(B)$. Moreover, $\lambda$ must be in one of the holes in $\sigma_{FLI}(A)$, otherwise equality (\ref{JEDNAKOSTSAETOM}) implies invertibility of $A-\lambda$, a contradiction. This finishes our proof. $\square$

\begin{Posledica}
    Let $A\in\mathcal{B}(\mathcal{H}),\ B\in\mathcal{B}(\mathcal{K})$. If $\sigma_d(A)\cap\sigma_l(B)$ has no interior points, then
$$\sigma_{FLI}(A)\cup\sigma_{FLI}(B)=\sigma_{FLI}(M_C)$$
    for every $C\in\mathcal{B}(\mathcal{K},\mathcal{H})$.
    Specially, if $A$ or $B$ is a compact operator, then equality holds for every $C\in\mathcal{B}(\mathcal{K},\mathcal{H})$.
\end{Posledica}
\textbf{Proof. }The first part of the proof follows directly from Theorem \ref{TEOREMAZAFLI}. The second part follows from the fact that a spectrum of a compact operator is at most countable. $\square$\\

With a little help of Theorem \ref{GLAVNA1}, we can even specify the precise form of the holes mentioned in Theorem \ref{TEOREMAZAFLI}. 

\begin{Teorema}\label{SESTAn=2}
Let $A\in\mathcal{B}(\mathcal{H}),\ B\in\mathcal{B}(\mathcal{K})$. Then
\begin{equation}\label{SESTAJEDAN}
\sigma_{FLI}(A)\cup\sigma_{FLI}(B)=\sigma_{FLI}(M_C)\cup W.
\end{equation}
holds for every $C\in\mathcal{B}(\mathcal{K},\mathcal{H})$, where
\begin{equation*}
\begin{aligned}
    W=\{\lambda\in\rho_{l}(A)\cap\rho_{SF-}(B):\ \alpha(B-\lambda)=\beta(A-\lambda)=\infty\}\cup\\ \lbrace\lambda\in\rho_{l}(A)\cap\rho_{SF-}(B):\ 0<\alpha(B-\lambda)\leq\beta(A-\lambda)<\infty\rbrace.
    \end{aligned}
    \end{equation*}
\end{Teorema}
\textbf{Proof.} It is not hard to see that the right side is a subset of the left side. Now, we prove the other inclusion. Put $W_1=\{\lambda\in\rho_{l}(A)\cap\rho_{SF-}(B):\ \alpha(B-\lambda)=\beta(A-\lambda)=\infty\}$ and $W_2=\lbrace\lambda\in\rho_{l}(A)\cap\rho_{SF-}(B):\ 0<\alpha(B-\lambda)\leq\beta(A-\lambda)<\infty\rbrace$, so that $W=W_1\cup W_2$. Assume that $\lambda\not\in\sigma_{FLI}(M_C)\cup W$. Since $\lambda\not\in\sigma_{FLI}(M_C)$, by Theorem \ref{GLAVNA1} we have that $\lambda\not\in\sigma_{FLI}(A)$, so it remains $\lambda\not\in\sigma_{FLI}(B)$ to be proven. By Theorem \ref{GLAVNA1} we also have that $\beta(B-\lambda)<\infty$ and that 
\begin{equation}\label{USLOVUTEOREMI}
\alpha(B-\lambda)\leq\beta(A-\lambda)<\infty\ or\ \alpha(B-\lambda)=\beta(A-\lambda)=\infty.
\end{equation} Since $\lambda\not\in W_1$ we conclude that either $\alpha(B-\lambda)\neq\beta(A-\lambda)$ or $\alpha(B-\lambda)=\beta(A-\lambda)<\infty$. Both of these cases exclude the second possibility in \eqref{USLOVUTEOREMI}, thus from \eqref{USLOVUTEOREMI} we have that 
\begin{equation}\label{USLOVUTEOREMI2}
\alpha(B-\lambda)\leq\beta(A-\lambda)<\infty.
\end{equation}
Since $\lambda\not\in W_2$, we conclude that $\alpha(B-\lambda)>\beta(A-\lambda)$ or $\alpha(B-\lambda)=0$. The first possibility is not possible because of \eqref{USLOVUTEOREMI2}, therefore $\alpha(B-\lambda)=0$ and so $\lambda\not\in\sigma_{FLI}(B)$. $\square$
\begin{Primedba}
    One can also write
    \begin{equation*}
\begin{aligned}
    W=\{\lambda\in\rho_{l}(A)\cap\rho_{SF-}(B):\ \alpha(B-\lambda)=\beta(A-\lambda)=\infty\}\cup\\ \lbrace\lambda\in\rho_{FLI}(A)\cap\rho_e(B):\ 0<\alpha(B-\lambda)\leq\beta(A-\lambda)\rbrace.
    \end{aligned}
    \end{equation*}
\end{Primedba}
\begin{Posledica}\cite{DONG2}\label{DONGOVATEOREMA}
Let $A\in\mathcal{B}(\mathcal{H}),\ B\in\mathcal{B}(\mathcal{K})$. Then 
$$
\sigma_{FLI}(A)\cup\sigma_{FLI}(B)=\sigma_{FLI}(M_C)
$$
for every $C\in\mathcal{B}(\mathcal{K},\mathcal{H})$ if
\begin{equation*}
\begin{aligned}
\lbrace\lambda\in\rho_{l}(A)\cap\rho_{SF-}(B):\ \alpha(B-\lambda)=\beta(A-\lambda)=\infty\}=\\ \lbrace\lambda\in\rho_{l}(A)\cap\rho_{SF-}(B):\ 0<\alpha(B-\lambda)\leq\beta(A-\lambda)<\infty\rbrace=\emptyset.
\end{aligned}
\end{equation*}
\end{Posledica}
Notice that we have obtained Theorem 3.10 from \cite{DONG2}.\\

By duality, we are now in a position to state analogous results concerning 'filling in holes' for the Fredholm right spectrum of $M_C$.

\begin{Lema}\label{IZMEDJUCEGAJESPEKTAR2}
    Let $A\in\mathcal{B}(\mathcal{H}),\ B\in\mathcal{B}(\mathcal{K})$. Then
    \begin{equation}
    \sigma_{FRI}(B)\setminus(\sigma_p(B)\cap\sigma_r(A))\subseteq\sigma_{FRI}(M_C)\subseteq\sigma_{FRI}(A)\cup\sigma_{FRI}(B)
    \end{equation}
    for every $C\in\mathcal{B}(\mathcal{K},\mathcal{H})$.
\end{Lema}
\textbf{Proof. }The second inclusion easily follows by \cite[Proposition 3.1]{OPERATORTHEORY} and factorization \eqref{RAZLAGANJE}. To prove the first inclusion, assume  $\lambda\not\in\sigma_{FRI}(M_C)$. Let $\lambda\in\sigma_{FRI}(B)\setminus(\sigma_p(B)\cap\sigma_r(A))$. Then $M_C-\lambda$ is a FRI operator, implying that $B-\lambda$ is right invertible. But, $\lambda\in\sigma_{FRI}(B)$, therefore $\alpha(B-\lambda)=\infty$. This condition, together with $\lambda\not\in\sigma_p(B)\cap\sigma_r(A)$ gives $\lambda\in\rho_r(A)$. However, this means that $M_C-\lambda$ and $A-\lambda$ are both Fredholm, implying that $B-\lambda$ must be Fredholm as well (see Lemma \ref{POMOCNALEMA}(e)). But, this is contrary to $\alpha(B-\lambda)=\infty$. Thus, the first inclusion is proved as well. $\square$ 

\begin{Teorema}\label{TEOREMAZAFRI}
    Let $A\in\mathcal{B}(\mathcal{H}),\ B\in\mathcal{B}(\mathcal{K})$. Then, for every $C\in\mathcal{B}(\mathcal{K},\mathcal{H})$,
    $$\eta(\sigma_{FRI}(A)\cup\sigma_{FRI}(B))=\eta(\sigma_{FRI}(M_C)),$$
    where $\eta(\cdot)$ denotes the 'polynomially-convex hull'.
    More precisely,
    $$\sigma_{FRI}(A)\cup\sigma_{FRI}(B)=\sigma_{FRI}(M_C)\cup W,$$
    where $W$ is equal to the union of some holes in $\sigma_{FRI}(B)$ which happen to be subsets of $\sigma_p(B)\cap\sigma_r(A)$.
\end{Teorema}

\textbf{Proof. } In the same way as in the proof of Theorem \ref{TEOREMAZAFLI} we have that
\begin{equation}\label{JEDNAKOSTSAETOM2}
    \eta(\sigma(T))=\eta(\sigma_{FRI}(T)).
\end{equation}
Assume $\lambda\in(\sigma_{FRI}(A)\cup\sigma_{FRI}(B))\setminus\sigma_{FRI}(M_C)$. Then, according to the first inclusion in Lemma \ref{IZMEDJUCEGAJESPEKTAR2} we have $\lambda\in\sigma_p(B)\cap\sigma_r(A)$. Moreover, $\lambda$ must be in one of the holes in $\sigma_{FRI}(B)$, otherwise equality (\ref{JEDNAKOSTSAETOM2}) implies invertibility of $B-\lambda$, a contradiction. This finishes our proof. $\square$

\begin{Posledica}
    Let $A\in\mathcal{B}(\mathcal{H}),\ B\in\mathcal{B}(\mathcal{K})$. If $\sigma_p(B)\cap\sigma_r(A)$ has no interior points, then
$$\sigma_{FRI}(A)\cup\sigma_{FRI}(B)=\sigma_{FRI}(M_C)$$
    for every $C\in\mathcal{B}(\mathcal{K},\mathcal{H})$.
    Specially, if $A$ or $B$ is a compact operator, then equality holds for every $C\in\mathcal{B}(\mathcal{K},\mathcal{H})$.
\end{Posledica}

In the following result we will specify the exact form of the holes mentioned in Theorem \ref{TEOREMAZAFRI}.

\begin{Teorema}
Let $A\in\mathcal{B}(\mathcal{H}),\ B\in\mathcal{B}(\mathcal{K})$. Then
\begin{equation}
\sigma_{FRI}(A)\cup\sigma_{FRI}(B)=\sigma_{FRI}(M_C)\cup W.
\end{equation}
holds for every $C\in\mathcal{B}(\mathcal{K},\mathcal{H})$, where
\begin{equation*}
\begin{aligned}
    W=\{\lambda\in\rho_{r}(B)\cap\rho_{SF+}(A):\ \alpha(B-\lambda)=\beta(A-\lambda)=\infty\}\cup\\ \lbrace\lambda\in\rho_{r}(B)\cap\rho_{SF+}(A):\ 0<\beta(A-\lambda)\leq\alpha(B-\lambda)<\infty\rbrace.
    \end{aligned}
    \end{equation*}
\end{Teorema}
\textbf{Proof. }This Theorem is proved using Theorem \ref{GLAVNA2}. The proof is very much similar to the proof of Theorem \ref{SESTAn=2}, and therefore we will omit it. $\square$\\

\begin{Primedba}
    One can also write
    \begin{equation*}
\begin{aligned}
    W=\{\lambda\in\rho_{r}(B)\cap\rho_{SF+}(A):\ \alpha(B-\lambda)=\beta(A-\lambda)=\infty\}\cup\\ \lbrace\lambda\in\rho_{FRI}(B)\cap\rho_e(A):\ 0<\beta(A-\lambda)\leq\alpha(B-\lambda)\rbrace.
    \end{aligned}
    \end{equation*}
\end{Primedba}
\begin{Posledica}
Let $A\in\mathcal{B}(\mathcal{H}),\ B\in\mathcal{B}(\mathcal{K})$. Then 
$$
\sigma_{FRI}(A)\cup\sigma_{FRI}(B)=\sigma_{FRI}(M_C)
$$
for every $C\in\mathcal{B}(\mathcal{K},\mathcal{H})$ if
\begin{equation*}
\begin{aligned}
\lbrace\lambda\in\rho_{r}(B)\cap\rho_{SF+}(A):\ \alpha(B-\lambda)=\beta(A-\lambda)=\infty\}=\\ \lbrace\lambda\in\rho_{r}(B)\cap\rho_{SF-}(A):\ 0<\beta(A-\lambda)\leq\alpha(B-\lambda)<\infty\rbrace=\emptyset.
\end{aligned}
\end{equation*}
\end{Posledica}

Finally, we point out some special classes of operators that enable us to simplify many of previously stated results. To this end consider classes
\begin{equation*}
\begin{aligned}
\mathcal{S}_+(\mathcal{H},\mathcal{K})=\lbrace T\in\mathcal{B}(\mathcal{H},\mathcal{K}):\ (\forall\lambda\in\mathds{C})\alpha(T-\lambda)\geq\beta(T-\lambda)\\ if\ at\ least\ one\ of\ these\ quantities\ is\ finite\rbrace,
\end{aligned}
\end{equation*}
\begin{equation*}
\begin{aligned}
\mathcal{S}_-(\mathcal{H},\mathcal{K})=\lbrace T\in\mathcal{B}(\mathcal{H},\mathcal{K}):\ (\forall\lambda\in\mathds{C})\alpha(T-\lambda)\leq\beta(T-\lambda)\\ if\ at\ least\ one\ of\ these\ quantities\ is\ finite\rbrace,
\end{aligned}
\end{equation*}
introduced in \cite{OPERATORTHEORY}. This classes are very important, in a sense they make many filling in holes results more simple (see e.g. \cite[Corollary 4]{HWANG},\cite[Theorem 5.7]{OPERATORTHEORY}). These classes are fairly wide and they contain many classes frequently seen in operator theory. For example, the class of quasihyponormal operators is contained in $\mathcal{S}_-$ \cite{OPERATORTHEORY}, while all quasitriangular (e. g. compact or cohyponormal operators) are part of $\mathcal{S}_+$ \cite{HWANG}.
\begin{Propozicija}
    Let $A\in\mathcal{S}_+(\mathcal{H})$, $B\in\mathcal{B}(\mathcal{K})$. Then 
$$
\sigma_{FLI}(A)\cup\sigma_{FLI}(B)=\sigma_{FLI}(M_C)
$$
for every $C\in\mathcal{B}(\mathcal{K},\mathcal{H})$.
\end{Propozicija}
\textbf{Proof. }This follows directly from Corollary \ref{DONGOVATEOREMA} and definition of $\mathcal{S}_+$. $\square$

Dually, one can obtain:

\begin{Propozicija}
    Let $A\in\mathcal{B}(\mathcal{H})$, $B\in\mathcal{S}_-(\mathcal{K})$. Then 
$$
\sigma_{FRI}(A)\cup\sigma_{FRI}(B)=\sigma_{FRI}(M_C)
$$
for every $C\in\mathcal{B}(\mathcal{K},\mathcal{H})$.
\end{Propozicija}

\textbf{Data availability statement}\\[3mm]
We declare there are no datasets associated with this work.\\

\textbf{Ethics declaration statement}\\[3mm]
Author declare there are no conflicts of interest associated with this work and that there is no relevant funding other than mentioned which influenced the writing of this manuscript.

\end{document}